\documentstyle[12pt]{article}

\begin{document}
\newcommand{\tr}{\bigtriangleup}
\newcommand{\epn}{\epsilon_q^{1\dots N}}
\newcommand{\bin}[2]{\left(\!\!\begin{array}{c}#1\\#2
\end{array}\!\!\right)}
\newcommand{\be}{\begin{equation}}
\newcommand{\ee}{\end{equation}}
\newcommand{\ba}{\begin{eqnarray}}
\newcommand{\ea}{\end{eqnarray}}
\def\lb{\label}
\def\theequation{\thesection.\arabic{equation}}

\begin{titlepage}

\title{ Characteristic Relations for Quantum Matrices }

\author{
P.~N.~Pyatov
\thanks{Supported in part by INTAS grant N\underline{o} 93-127.}
\thanks{e-mail address: pyatov@thsun1.jinr.dubna.su}
\\
\it Bogoliubov Laboratory of Theoretical Physics, \\
\it Joint Institute for Nuclear Research, \\
\it 141980 Dubna, Moscow Region, Russia \\
\rm and \\
P.~A.~Saponov
\thanks{e-mail address: saponov@mx.ihep.su}
\\
\it Theoretical Department, \\
\it Institute for High Energy Physics, \\
\it 142284 Protvino, Moscow Region, Russia
}
\date{}
\maketitle

\end{titlepage}

\newpage

\vspace*{5cm}
\begin{center}
\bf Abstract
\end{center}

General algebraic properties of the algebras of vector fields
over quantum linear groups $GL_q(N)$ and $SL_q(N)$ are studied.
These quantum algebras appears to be quite similar to the classical
matrix algebra. In particular, quantum analogues of the
characteristic polynomial and characteristic identity are obtained
for them. The $q$-analogues of the Newton relations connecting two
different generating sets of central elements of these algebras
(the determinant-like and the trace-like ones) are derived.
This allows one to express the $q$-determinant of quantized vector
fields in terms of their $q$-traces.

\newpage

\section{Introduction}
\setcounter{equation}0

Since the time of their discovery quantum groups were presented
in several closely related but not strictly equivalent forms.
Being originally obtained as the quantized universal enveloping
(QUE) algebras \cite{Dri,Jimbo} they were then reformulated
in a matrix form \cite{FRT}. In this latter approach a quantum
group is generated by a pair of upper- and lower-triangular
matrix generators $L_+$ and $L_-$ satisfying quadratic permutation
relations. A further variation of this approach is to combine
$L_+$ and $L_-$ into a single matrix generator $L = S(L_-)L_+$.
Here $S( \cdot )$ is the usual notation for the antipodal mapping.
Following Ref.\cite{Kul-Skl} we will call the algebra generated by
the matrix generator $L$ as the reflection equation algebra (REA).
After suitable completing this algebra can be related to the quantum
group by the (Hopf algebra) isomorphism \cite{Bur,DJSWZ}, although
the Hopf structure is implicit in REA formulation. This algebra have
found several applications (see \cite{Kul-Skl,Resh-ST,Kulish,Majid}
and references therein). Let us mention here only one of them.
Namely, in the construction of the quantum group differential
calculus the matrix generators $L$  are used as the basic set of
(right-)invariant vector fields \cite{Jurco,Aleks-Fad,Zum}.

A remarkable property of the REA formulation is that the algebra of
the $L$ matrices in several aspects appears to be quite similar to
the classical matrix algebra. In particular, both the notions of the
matrix trace and that of the matrix determinant admit the
generalization to the case of $L$ matrices (see
\cite{FRT,Resh} and \cite{DJSWZ,SWZ,FP}). In the present paper we
intend to establish further similarities of the REA with the
classical matrix algebra. We restrict ourselves to considering
the REA of the  $GL_q(N)$ and/or $SL_q(N)$ type only. For these
algebras the recurrent formulas relating two different center
generating sets, the determinant-like and the trace-like ones,
are obtained. These formulas are the quantum analogues of the
classical Newton relations. Further we define the characteristic
polynomial and derive the characteristic identities (the analogue
of the Cayley-Hamilton theorem) for the $L$-matrices. Existence of
such identities were first mentioned  in Ref. \cite{Kul-Skl} where
they were presented for the case $N = 2$.
For general $N$ the characteristic identities were obtained
in Ref.\cite{Naz-Tar}
(see Remark 4.8  of this paper)
when studying the algebraic structure
of quantum Yangians $Y_q({\sl gl}(N))$. We reproduce this result
in REA approach. Then, by the joint use of the quantum
Newton relations and the characteristic identities
one can obtain the expressions for
trace-like central elements of higher powers.
In the QUE presentation
the similar characteristic identities were considered in
Ref.\cite{GZB}.
We believe that the REA presentation and the use of
$R$-matrix technique makes all considerations and the final formulas
much clearer.

We conclude this Section with a brief remainder on some facts from
the classical theory of matrices (see e.g. \cite{Lank}) which then
will be generalized to the quantum case.

Consider $N\times N$ matrix $A$ with complex entries. Its
characteristic polynomial is defined as
\be
\lb{char}
\tr(\,x\,) \equiv \det ||\, x\,{\bf 1} - A\,||
\equiv x^N + \sum_{k=1}^{N} (-1)^k \sigma(k) \,x^{N-k} \; .
\ee
The eigenvalues $\{\lambda_i\}$, $i=1, \ldots ,N$, of the matrix $A$
are  solutions of the characteristic equation $\tr(x)=0$. The
coefficients $\sigma(i)$ of the characteristic polynomial if expressed
in terms of $\lambda_i$ form the set of basic symmetric polynomials
of $N$ variables:
\be
\lb{sigmas}
\sigma(1) \equiv \sum_{i=1}^{N} \lambda_i  =  Tr\,A \; ,\;\;\;\;
\sigma(2) \equiv \sum_{i<j} \lambda_i \lambda_j \;\; , \, \ldots \, ,
\;\;\;\;
\sigma(N) \equiv \prod_{i=1}^{N} \lambda_i = \det A \; .
\ee
Also one can directly express $\sigma(i)$ in terms of the matrix
elements of $A$. Up to some numerical factor each $\sigma(i)$ is
given by the sum of all the principal minors of the $i$-th order
\be
\lb{sigmas2}
\sigma(i) = { 1 \over i! \, (N-i)!  } \,
\epsilon^{1\dots N} A_1\dots A_i\,\epsilon^{1\dots N} \;.
\ee
Here $\epsilon^{1\dots N}$ is the antisymmetric Levi-Civita
$N$-tensor. The compressed matrix conventions used in this
formula will be explained later (directly in the quantum case).

Another standard set of symmetric polynomials is given by the
traces of powers of the matrix $A$
\be
\lb{traces}
s(i) \equiv \sum_{k=1}^{N} (\lambda_k)^i =  Tr\,A^i \; , \;\;
1\le i \le N \; .
\ee
The two basic sets $\{\sigma(i)\}$ and $\{s(i)\}$ are connected by
the so-called Newton relations
\be
\lb{newton}
i\,\sigma(i) - s(1)\,\sigma(i-1) + \dots + (-1)^{i-1} s(i-1) \,
\sigma(1) + (-1)^i s(i) = 0.
\ee
In particular these recurrent relations allow one  to express the
determinant of the matrix $A$ as a polynomial of its powers traces.

Finally, if one substitutes the matrix $A$ into the characteristic
polynomial (\ref{char}) instead of the scalar variable $x$  then
the resulting matrix expression vanishes identically. This is the
Cayley-Hamilton theorem, and according to it any function of matrix
$A$ can be reduced to a polynomial of an order not exceeding $N-1$.

\section{Quantum Newton Relations and Characteristic Polynomial}
\setcounter{equation}0

First of all let us introduce some definitions and notation to be
used in what follows. The REA is defined as the algebra generated by
matrix generators $L$ subject to the following permutation rules
\be
\lb{rlrl}
L_1 \, {\hat R}_{12} \, L_1 \, {\hat R}_{12} \, = \,
{\hat R}_{12} \, L_1 \, {\hat R}_{12} \, L_1 \; .
\ee
Here the standard conventions for denoting matrix spaces (see
\cite{FRT}) are used. $\hat{R}_{12}$ is the $GL_q(N)$ $R$-matrix
\cite{Jimbo} satisfying the Yang-Baxter equation and the Hecke
condition respectively
\ba
\lb{YB}
{\hat R}_{12} {\hat R}_{23} {\hat R}_{12} &=&
{\hat R}_{23} {\hat R}_{12} {\hat R}_{23} \; , \\
\lb{H}
{\hat R}^2 &=& {\bf 1} + \lambda {\hat R} \; ,
\ea
where ${\bf 1}$ is the unit matrix, and $\lambda = q - 1/q$. Below
we will further compress this notations denoting ${\hat R}_{i (i+1)}
\equiv {\hat R}_i$ and omitting the index of $L$-matrix: $L_1
\equiv L$,  since it always appears in the first matrix space.

We will also need the notions of the quantum trace and the
$q$-deformed Levi-Civita tensor. The operation $Tr_q$ of taking
the quantum trace of $N \times N$ quantum matrix $X$ looks like
\be
\lb{tr}
Tr_q(X) = Tr({\cal D}X) \;\;\; , \;\;\;
{\cal D} = diag \{ q^{-N+1}, q^{-N+3}, \ldots ,  q^{N-1} \} \; .
\ee
The $q$-deformed Levi-Civita tensor
$\epsilon_{q}^{i_1 \ldots i_N}$ (or $\epsilon_{q}^{ 1 \ldots N}$
in brief notation) is defined, up to a factor by its characteristic
property
\be
\lb{e}
\left( {\hat R}_i + {1 \over q} \right)
\epsilon_{q}^{1 \ldots N} = 0
\; , \;\;\;\; 1 \leq i \leq N-1 \; .
\ee
The normalization is usually fixed by demanding
$\epsilon_{q}^{i_1 \ldots i_N} = 1$ for $i_1 = 1,\ldots ,i_N = N$.
Its square is then equal to
$$
| \epsilon_q {|}^2 \equiv
\epsilon_{q}^{1 \ldots N} \epsilon_{q}^{1 \ldots N} =
q^{N(N-1)/2} {N_q}!,
$$
where $p_q = (q^p - q^{-p})/\lambda$ are the usual $q$-numbers.

In ref. \cite{FRT} the two generating  sets for the center of the
REA were presented. One of them is formed by the trace-like elements
\be
\lb{s}
s_q(i) = q^{1-N} Tr_q L^i \qquad 1 \leq i \leq N\, ,
\ee
where the normalizing factor is chosen for future convenience.
Another generating set consists  of the determinant-like  elements
$
 \sim \epsilon_q^{1 \ldots N} {L_-}_N \dots {L_-}_{(i+1)}{L_+}_i,
\dots {L_+}_1\epsilon_q^{1 \ldots N} .
$
For our purposes it is better to express these
generators in terms of the $L$-matrices:
\be
\sigma_q(i) = \alpha_i\epn (L_1 {\hat R}_1 \dots
{\hat R}_{i-1})^i\epn. \label{sig}
\ee
Here $\alpha_i$ are some  normalizing constants. First of them
is fixed to be
$$
\alpha_1 = q^{1-N} N_q / |\epsilon_q|^2
$$
by the natural condition $\sigma_q(1) = s_q(1)$. The other ones
will be specified below.

The connection between the two basic sets \{$\sigma_q(k)$\} and
\{$s_q(k)$\} is provided by quantum analogues of the Newton relations
(\ref{newton}). At the classical level their derivation is usually
based on using the spectral presentations (\ref{sigmas}),
(\ref{traces}) for $\sigma(k)$ and $s(k)$. However this presentation
is not available in quantum case. Indeed, the spectrum of the quantum
matrix $L$ can be constructed only if the center of REA is
algebraically closed. The latter in turn can be treated only on
passing to a concrete representation of REA. In order to by-pass
this difficulty we shall develop some more technique.

Define an operator $S_N$ which symmetrizes any $N\times N$ matrix
$X$ in $N$ matrix spaces:
\be
S_N(X) = X_1 + {\hat R}_1 X_1 {\hat R}_1 + \dots
+ {\hat R}_{N-1}\dots {\hat R}_1 X_1 {\hat R}_1\dots {\hat R}_{N-1}.
\label{sn}
\ee
The characteristic properties of this symmetrizer
\be
\left[S_N(X)\,,\,{\hat R}_i\right] = 0\, , \;\;\;
1 \leq i \leq N-1 \; , \label{chpr}
\ee
are fulfilled due to the relations (\ref{YB}), (\ref{H}).
Further, the following useful formula
\be
\label{as}
\epsilon_q^{1\dots N} S_N(X) =
s_X \epsilon_q^{1\dots N}
\ee
is a direct consequence of (\ref{chpr}) and (\ref{e}). Here the
scalar factor $s_X$ is
$$
s_X = \frac{1}{|\epsilon_q|^2}\epn S_N(X)
\epn = q^{1-N} Tr_qX.
$$
For $X = L^i$ this factor coincides with the $s_q(i)$ (\ref{s}).
Now we are able to  prove \medskip

{\bf Proposition.}\hspace*{5 true mm}
\it
For $1 \leq i \leq N$
the generators
 $\sigma_q(i)$ and $s_q(i)$ are connected by the relations
\be
\lb{qnewton}
{i_q \over q^{i-1}} \sigma_q(i) -
s_q(1)\,\sigma_q(i-1) + \dots +
(-1)^{i-1} s_q(i-1) \,\sigma_q(1) +
(-1)^i s_q(i) = 0 \, ,
\ee
provided that the numerical factors $\alpha_i$ are fixed as follows
\rm
\be
\label{alpha}
\alpha_i = {N_q! \over (N-i)_q!\,i_q!} \;
{q^{-i(N-i)} \over |\epsilon_q|^2 }  \; .
\ee
{\bf Proof.}
Consider the quantities $s_q(i-p)\sigma_q(p)$ for $1 \leq p \leq i-1$.
With the help of (\ref{as}) and
the definitions of $s_q(i)$ and $\sigma_q(i)$
one can perform the following transformations:
\begin{eqnarray*}
s_q(i-1) \sigma_q(1) & = & \alpha_1 \underline{s_q(i-1) \epn}
L\,\epn = \alpha_1 \epn S_N(L^{i-1}) L\, \epn \\
& = & s_q(i) + \alpha_1 {(N-1)_q \over q^{N-2}} \epn
L^{i-1} R_1 L\, R_1 \epn \; ;
\\
& & ................................
\\
s_q(i-p) \sigma_q(p) & = &
\alpha_p {p_q \over q^{p-1}} \epn (L^{i-p+1} R_1\dots R_{p-1})
(L\,R_1 \dots R_{p-1})^{p-1} \epn  \\
& + &
\alpha_p {(N-p)_q \over q^{N-p-1}} \epn (L^{i-p}R_1\dots R_p)
(L\,R_1\dots R_p)^p \epn  \; ;
\\
& & ................................
\\
s_q(1) \sigma_q(i-1) & = &
\alpha_{i-1} {(i-1)_q \over q^{i-2}} \epn (L^{2} R_1\dots R_{i-2})
(L\,R_1 \dots R_{i-2})^{i-2} \epn  \\
& + &
\alpha_{i-1} {(N-i+1)_q \over q^{N-i}} \;
{\sigma_q(i) \over \alpha_i} \; .
\end{eqnarray*}
Now the arbitrary coefficients $\alpha_p$ should be fixed in such a
way that the last term in $s_q(i-p+1)\sigma_q(p-1)$  and the first
one in $s_q(i-p)\sigma_q(p)$  be equal. This is the case if $\alpha_p$
satisfy the relations
\be
\lb{recurr}
\alpha_p = q^{2p-1-N} { (N-p+1)_q \over p_q } \alpha_{p-1}\; .
\ee
Then, on taking the alternating sum $\sum_{p}(-1)^{p-1}s_q(i-p)
\sigma_q(p)$ we find that the only terms which survive are the first
one in $s_q(i-1)\sigma_q(1)$ and the last one in $s_q(1)\sigma_q(i-1)$
and, thus, we obtain relations (\ref{qnewton}). Finally, given the
value of $\alpha_1$ one easily shows that recursion (\ref{recurr}) is
solved by (\ref{alpha}).
{}~~\rule{2.5mm}{2.5mm}\par
Few remarks are in order here:\par
{\bf 1.} Up to now we always consider the $GL_q(N)$ matrices.
Specification to the $SL_q(N)$ case can be achieved by fixing
quantum determinant of the $L$-matrix (see \cite{DJSWZ,SWZ,FP}):
$Det\, L = q^{1-N} \sigma_q(N) = 1$.\par
{\bf 2.} It is worth mentioning that the singular points where
the connection between $\{s_q(i)\}$ and $\{\sigma_q(i)\}$
breaks down are the roots of unity: $k_q = 0$ for $1 \leq k \leq N$.
This is apparently related to the fact that the isomorphism of
the Hecke algebra of $A_{N-1}$ type and the group algebra of
symmetric group  ${\bf C}S_N$ is also destroyed at these points
(see e.g. \cite{Wenzl}).
\medskip

Now let us turn to deriving the quantum characteristic identity
for the matrix $L$. It can be found in a way quite similar to
that of the classical case. Namely, we should find a matrix
polynomial $B$ of $(N-1)$-th order of $L$ obeying the relation
\be
\lb{B}
(L - x{\bf 1})\, B(L,x)\, \epn = \epn \tr(x) \; .
\ee
Here $x$ is a ${\bf C}$-number variable and $\tr(x)$ is the
scalar polynomial of $x$  -- the characteristic polynomial of the
$L$-matrix. The following is a generalization of the
Cayley-Hamilton theorem to the quantum case.\par
{\bf Theorem.}\hspace*{5mm}
\it
The matrix polynomial $B(L,x)$ when defined as
\be
\lb{defB}
B(L,x) = {\hat R}_1\dots {\hat R}_{N-1} \prod_{i=1}^{N-1}
\left[ (L-q^{2i}x{\bf 1}){\hat R}_1\dots {\hat R}_{N-1} \right] \; ,
\ee
satisfies the relation {\rm (\ref{B})}. The characteristic polynomial
of the matrix $L$ looks like
\be
\lb{polynom}
\tr(x) = \sum^{N}_{i=0} (-x)^i \sigma_q(N-i) \; ,
\ee
and for the $L$-matrix the following characteristic
identity is satisfied
\be
\tr(L) = \sum^{N}_{i=0} (-L)^i \sigma_q(N-i) \equiv 0 \; .
\ee
\rm
\par
{\bf Proof:} The relation (\ref{B}) is fulfilled if and only if
its left hand side is totally $q$-antisymmetric, i.e. if it obeys
the characteristic relations (\ref{e}) of the $q$-antisymmetric
tensor. This in turn is valid if the matrix quantity
$(L - x{\bf 1})\, B(L,x)$ commutes with all the ${\hat R}_i$,
$1 \leq i \leq N-1$ up to terms proportional to the
$q$-symmetric projectors ${P_+}_i \equiv ( {\hat R}_i +1/q)/2_q$,
which vanish when contracting with $\epn$.
Now the key observation is that the commutator
\mbox{$[{\hat R}_1\,,\,(L-x{\bf 1}){\hat R}_1(L-\beta x{\bf 1})
{\hat R}_1 ]$} is proportional to $P_+$ only if $\beta = q^2$.
With this observation the construction of $B$ becomes clear
and one immediately checks that chosen as in (\ref{defB}) $B$
does satisfy the relation (\ref{B}).

Then as a direct consequence of (\ref{B}) and (\ref{defB})
we get the following expression for $\tr(x)$:
$$
\tr(x) = {1 \over |\epsilon|^2}
\epn \prod_{i=0}^{N-1} \left[ (L -q^{2i} x {\bf 1})
{\hat R}_1 \dots {\hat R}_{N-1} \right] \epn \; .
$$
This expression can be further simplified with the use of
(\ref{e}), (\ref{rlrl}), (\ref{YB}) together with the
$q$-combinatorial relations. The calculations are straightforward but
rather lengthy and we omit them here presenting
the result in (\ref{polynom}).

To prove the characteristic identity we shall contract
relation (\ref{B}) with $\epsilon_q^{2\dots N+1}$:
$$
(L - x{\bf 1})\,
\epsilon_q^{2\dots N+1} B(L,x)\, \epn =
(\epsilon_q^{2\dots N+1} \epn) \tr(x) \; .
$$
The right hand side of this relation is proportional to the
unit matrix and, hence, $\epsilon_q^{2\dots N+1} B \epn$
is proportional to $(L - x{\bf 1})^{-1}$. The classical limit of
this relation is the standard base for proving the
Cayley-Hamilton theorem \cite{Lank}. In quantum case all the
considerations are completely the same, and the resulting statement
is that the matrix polynomial $\tr(L)$ identically vanishes.
\rule{2.5mm}{2.5mm}

Here are few final comments.\par
{\bf 1.} The characteristic identity provides us with the
compact expression for the inverse matrix of $L$:
$$
L^{-1} = {1 \over \sigma_q(N)}
\sum_{i=0}^{N-1} (-L)^{i} \sigma_q(N-i-1) \; .
$$

{\bf 2.} Multiplying the characteristic identity by $L^p$,
and taking the $q$-trace one gets the expressions of higher
symmetric polynomials $s_q(N+p)$
in terms of the basic ones $\sigma_q(i)$.\par
{\bf 3.} On passing to concrete REA representations
the order of the characteristic identity may decrease
due to the basic symmetric
polynomials $\sigma_q(i)$ become dependent.
This is illustrated in a recent paper \cite{Chu-Zum}
where the $L$-matrices were realized as pseudo-differential
operators acting on the quantum plane and they were found to
possess the characteristic identity of the second order.

\end{document}